\documentclass[10pt]{article}

\usepackage[centertags]{amsmath}
\usepackage{amsfonts,amssymb,amsthm}
\usepackage{authblk}
\usepackage{comment}
\usepackage{caption,graphicx}
\usepackage[text={6.0in,8.6in},centering,letterpaper]{geometry}
\usepackage{hyperref}
\usepackage{mathrsfs}
\usepackage{mathtools}
\usepackage[numbers,comma,square,sort&compress]{natbib}
\usepackage{subcaption}
\usepackage{tikz}
\usepackage{url}

\newcommand{\U}{\ensuremath{\Upsilon}}
\newcommand{\s}{\ensuremath{\sigma}}
\newtheorem{remark}{Remark}
\numberwithin{equation}{section}
\newcommand{\sref}[1]{(\ref{#1})}                       

\begin{document}
\author{C.H.S. Hamster}
\affil{Dutch Institute for Emergent Phenomena,\\University of Amsterdam.\\
Email:  {\normalfont{\texttt{c.h.s.hamster@uva.nl}}}}
\title{Right and Wrong Ans\"atze for Nonlinear Waves in Stochastic PDEs}
\date{\today}
\maketitle

\begin{abstract}
I investigate the possibility that explicit solutions of stochastic reaction-diffusion equations can be found by multiplying the deterministic travelling waves with a stochastic exponent. This approach has become widespread in the literature in recent years. I will conclude that this approach is, in general, not a valid Ansatz and only works in the case of NLS-type equations in the Stratonovich interpretation. 
\end{abstract}

\section{Introduction}
In the recent decade, efforts to study nonlinear waves in a stochastic setting have picked up pace, see e.g. \cite{kuehnreview, Hamster2020,cartwright2019,westdorp2024long,vanWinden2024noncommutative,bosch2024multidimensional, lang, Inglis, lord2012}. In the wake of these efforts, an industry evolved around computing explicit solutions to nonlinear SPDEs where the deterministic version is known to have a travelling wave, based on the following assumption: A travelling wave solution $u$ to the SPDE can be found by multiplying the deterministic travelling wave $(\Phi,c)$ with an appropriately chosen exponent depending on the noise and time. For an SPDE forced by a Brownian motion $\beta_t$ with intensity $\sigma$, this Ansatz could be
\begin{align}
    u(x,t)=\Phi(x-ct)e^{\sigma\beta_t-\frac{\s^2}{2}t}.
\end{align}
I will refer to this as Ansatz \textbf{A}, and we will study it in more detail later. It should be immediately clear that this approach cannot work in general, as the effect of the noise can, again generally, not be captured in a single 1D function.

To highlight this point, I will briefly discuss two stochastic travelling waves that can be explicitly computed. First, a moving front in a stochastic Nagumo equation:
\begin{align}
\label{eq:int:StochNag}
    du=[u_{xx}+u(1-u)(u-a)]dt+\sigma u(1-u)d\beta_t.
\end{align}
Here, $\beta_t$ is a 1D Brownian motion and $0<a<1$. It has a deterministic travelling wave 
\begin{align}
    u(x,t)=\frac{1}{2}\left[1+\tanh\left(-\frac{\sqrt{2}}{4}(x-ct)\right)\right], \hspace{5mm} c=\sqrt{2}\left(\frac{1}{2}-a\right).
\end{align}
Traditionally, this solution is found using a Travelling Wave Ansatz, i.e. one assumes that a fixed profile $\Phi$ exists that moves with speed $c$. Upon introducing the travelling wave coordinate $\xi=x-ct$, the pair $(\Phi,c)$ solves the Travelling Wave Equation (TWE)
\begin{align}
    \Phi''+c\Phi'+\Phi(1-\Phi)(\Phi-a)=0,
\end{align}
with appropriate boundary conditions. For the stochastic equation \sref{eq:int:StochNag}, an Ansatz is introduced in \cite{hamster2017}, where the assumption that $\Phi$ moves with speed $c$ is replaced with the assumption that a pair $(\Phi_\s,c_\s)$ exists such that $\Phi_\s(x-c_\s t-\alpha_t)$ is a solution for some stochastic process $\alpha_t$. This Ansatz results in the following solution:
\begin{align}
\begin{split}
    \Phi_\s(\xi)=&\Phi\left(\sqrt{1+\sigma^2}\xi\right), \hspace{3mm}
    c_\s=\frac{c}{\sqrt{1+\s^2}}, \hspace{3mm}   \alpha_t=\frac{\sigma\sqrt{2}}{\sqrt{1+\s^2}}\beta_t.
\end{split}
\end{align}
Note that this explicit solution can only be found when the noise is 1D and the nonlinear term is given by $g(u)\sim u(1-u)$. This quadratic function works for a very specific reason, as $\Phi(1-\Phi)\sim \Phi'$. This allowed the authors to construct a solution where $\Phi_\s(1-\Phi_\s)\sim\Phi'_\s$, which not only significantly reduces the complexity of the equation, but also ensures that the noise $g(\Phi_\s)d\beta_t$ only acts along the manifold of translates of $\Phi_\s$.
To my knowledge, there is no function $f(t,\beta_t)$ that would allow us to rewrite
\begin{align}
    \Phi_\s(x-c_\s t-\alpha_t)=\tilde\Phi_\s(x-\tilde c_\s t)f(t,\beta_t), 
\end{align}
for some deterministic, time-independent, function $\tilde\Phi_\s$, as Ansatz \textbf{A} requires.

A second example where an explicit stochastic travelling wave is known is the Korteweg-De Vries equation: 
\begin{align}
    du=[6uu_x-u_{xxx}]dt+\sigma d\beta_t.
\end{align}
Again, we assume that $\beta_t$ is a 1D Brownian motion. For this equation, an explicit solution is known \cite{wadati1983stochastic}:
\begin{align}
    u(x,t)&=\Phi_c\left(x-ct+6\sigma\int_0^t\beta_sds\right)+\sigma\beta_t,
\end{align}
where $\Phi_c$ is the solution to the deterministic TWE
\begin{align}
    \Phi_c'''-6\Phi_c\Phi_c'-c\Phi_c'=0, \hspace{3mm} \Phi_c(-\infty)=\Phi_c(\infty)=0.
\end{align}
Again, the solution $u$ cannot be written, to the best of my knowledge, in the form of Ansatz \textbf{A}.
In \cite{cartwright2021collective}, this solution is recovered using a so-called collective coordinate approach. In the collective coordinate approach, it is essential that the position is assumed to be stochastic. It must have been, since this approach was inspired by the explicit solution above.  

Both exact stochastic solutions have in common that, yes, the typical shape of the pattern is retained, the $\tanh$ for the Nagumo equation and the $\text{sech}^2$ for the KdV equation, but the scale parameters of the solution, like position, width and amplitude, are stochastically varying, and these effect are not captured by a single function. In this short note, we will start with some examples of explicit solutions. Then follows a case study of \cite{baber2024exact}, followed by several types of classic PDEs and discuss if and when a simple Ansatz like Ansatz \textbf{A} could work for a stochastic version of the PDE. We will conclude that it only works for stochastic Nonlinear Schr\"odinger equations in Stratonovich interpretation, like Eq.~\ref{eq:NLSstrat}.

\subsection{Origins of flawed Ans\"atze}
It is difficult to exactly pinpoint the origins of these papers. The first time I encountered it was in \cite{mohammed2021exactSchrodinger}. Well, not that version, but a version that I reviewed and advised to reject in November 2020 for another journal. I wrote a very extensive report, not only outlining the fundamental flaws (It\^o vs Stratonovich), but also pointed out many small (typographical) mistakes. However, all these mistakes ended up in the published version \cite{mohammed2021exactSchrodinger}. This means that \cite{albosaily2020exact} is the first to appear, but closely followed by \cite{mohammed2021exactKonno}, \cite{mohammed2021exact} and \cite{mohammed2021exactSchrodinger}. These papers have in common that W.W. Mohammed is one of the authors.   

At the same university as W.W. Mohammed, Mansoura University in Egypt, we also find M.A.E. Abdelrahman and M.A. Sohaly. From 2017, they wrote a series of papers where they look at solutions to PDEs with random coefficients~\cite{abdelrahman2017solitary,abdelrahman2018development}. First, they `solve' the PDE, and then study a stochastic version simply by replacing a parameter with a random variable. It is not clear what the point of this approach is.  

\section{A systematic approach to exact solutions}
We can systematically look for exact solutions by following an approach that I will call `reverse phase tracking'. Reverse in the sense that we do not find an Ansatz for an SPDE, but rather start with an Ansatz and find an SPDE. 

\subsection{Advective 1D noise}
Suppose we have the following very general SPDE
\begin{align}
    \label{eq:ExactSolGeneralS}
    dU=[F(U,U_x,U_{xx},...)]dt+\sigma U_x\circ d\beta_t,
\end{align}
and introduce the phase
\begin{align}
    d\Gamma(t)=cdt+\sigma d\beta_t.
\end{align}
Then, $v(t)=U(\cdot-\Gamma(t),t)$ solves the PDE
\begin{align}
    v_t=F(v,v_x,v_{xx},..)+cv_x.
\end{align}
This follows directly from the classic chain rule because the noise is interpreted in the Stratonovich sense. Therefore, if the deterministic version of Eq.~\sref{eq:ExactSolGeneralS} has a travelling wave, we can immediately solve the stochastic version. When we look at the It\^o interpretation, we introduce a phase 
\begin{align}
    d\Gamma(t)=c_\s dt+\sigma d\beta_t
\end{align}
and find that $U(\cdot-\Gamma(t),t)$ solves the PDE
\begin{align}
    v_t=F(v,v_x,v_{xx},..)+c_\sigma v_x-\frac{\s^2}{2}v_{xx}.
\end{align}
Hence, we can exactly solve the It\^o version of the SPDE when we can solve the perturbed TWE equation above. For example, we can solve the following stochastic Nagumo equation:
\begin{align}
    dU=[dU_{xx}+U(1-U)(U-a)]dt+\sigma U_xd\beta_t,
\end{align}
because we now know that $\Phi_\s(\cdot-\Gamma(t),t)$ is a solution when
\begin{align}
    (d-\frac{1}{2}\sigma^2)\Phi''_{\s}+c_\sigma\Phi'_\s+\Phi_\s(1-\Phi_\s)(\Phi_\s-a)=0.
\end{align}
Solving this equation is straightforward when $d-\frac{1}{2}\s^2>0$. Note that the boundary conditions do not change, as $U_x$ vanishes at the deterministic endpoints 0 and 1. However, if we were to use this approach for the KdV equation, as was done in \cite{cartwright2021collective}, we solely get a negative diffusive term. We cannot solve the resulting TWE exactly\footnote{In \cite{obeidat2025exploration}, the authors claim that they found a solution for this perturbed TWE, but none of the solutions presented in that paper are solitons.}, but we still conclude that the soliton is unstable due to the negative diffusive term. 

In the case of the following stochastic Burgers-KdV equation
\begin{align}
    dU=[\delta U_{xxx}+\mu U_{xx}+\beta UU_x]dt+\sigma U_xd\beta_t,
\end{align}
we can solve this SPDE by solving the following perturbed TWE:
\begin{align}
    \delta \Phi'''_\s+(\mu-\frac{\s^2}{2}) \Phi''_\s+c_\s\Phi'_\s+\delta \Phi_\s\Phi_\s'=0.
\end{align}
This time, the equation is stable and solvable as long as $(\mu-\frac{\s^2}{2})>0$. 
Note that for this specific equation, the speed $c$ is uniquely determined by the endpoints of the wave $\Phi_\s$:
\begin{align}
    c_\s=-\frac{\beta}{2}(\Phi_\s(-\infty)+\Phi_\s(\infty)).
\end{align}
Hence, $c_\s=c_0$ when we keep the endpoind fixed, but the profile $\Phi_\s$ is changed. However, in \cite{adjibi2024exact}, this not accounted for. There, the distinction between $\mu$ and $\mu_{\rm eff}=\mu-\s^2/2$ is not made, and the boundary conditions are not fixed. For the deterministic travelling wave presented in \cite{adjibi2024exact}, the boundary conditions depend on $\mu$. Hence, when we want to understand how the noise changes the solutions without changing the boundary conditions, I'm not sure we can find an exact solution. 

More generally, suppose $U(x,t)$ is a smooth function of space and time, and $\Gamma(t)$ is a stochastic process given by
\begin{align}
    \Gamma(t)=\Gamma_0+\int_0^t\alpha(s)ds+\int_0^t\sigma(s)d\beta_s
\end{align}
and 
\begin{align}
    dU=[F(t,U,U_x,U_{xx},...)]dt+\sigma(t) U_x d\beta_t
\end{align}
Then, $U(\cdot-\Gamma(t),t)$ solves the PDE
\begin{align}
    v_t=F(t,v,v_x,v_{xx},..)+\alpha(t)v_x-\frac{\s(t)^2}{2}v_{xx}.
\end{align}
Hence, whether or not this is useful depends purely on luck or design. In \cite{flores2020exact}, it was shown that 
\begin{align}
    dU=[(\exp(t)+\frac{\sigma^2}{2})U_{xx}+\exp(t)UU_x]dt+\sigma U_xdW_t
\end{align}
can be solved, because the PDE
\begin{align}
    v_t=\exp(t)v_{xx}+\exp(t)vv_x
\end{align}
can be solved. However, it is not clear why we would be interested in either equation. Hence, this approach feels more like a hammer looking for a nail than a nail looking for a hammer.

\subsection{Additive 1D noise}
Let us study the following Burgers-KdV equation:
\begin{align}
    dU=[\delta U_{xxx}+\mu U_{xx}+\beta UU_x+\alpha U_x]dt+\sigma(t)dW_t.
\end{align}
Inspired by the exact solution of the KdV equation as studied in the introduction and~\cite{adjibi2024exact}, we choose 
\begin{align}
    U(t)=\Phi_\s\left(\cdot-c_\sigma t+\beta\int_0^t\int_0^s\sigma(s')ds'dW_s\right)+\int_0^t\sigma(s)dW_s.
\end{align}
This results in
\begin{align}
    \left[-c_\s+\beta\int_0^t\sigma(s)ds\right]\Phi'_\s=\delta\Phi'''_{\s}+\mu \Phi''_{\s}+\beta \left(\Phi_\s+\int_0^t\sigma(s)d\beta_t\right)\Phi'_\s,
\end{align}
which reduces to
\begin{align}
    -c_\s\Phi'_\s=\delta\Phi'''_{\s}+\mu\Phi''_{\s}+\beta\Phi_\s\Phi'_\s.
\end{align}
Hence, the TWE is independent of $\sigma(t)$, and the shifted deterministic wave $(\Phi_0,c_0)$ plus Brownian motion solves the SPDE. Note that this method works because of the structure of the PDE. When we would study a modified KdV-Burgers equation, e.g. by replacing $uu_x$ with $u^2u_x$, this Ansatz wouldn't work anymore, but you can add any type of nonlinearity depending on spatial derivatives.
\section{A case study}
In \cite{baber2024exact}, the authors\footnote{In private communications, the authors acknowledge that their analysis might be wrong, but chose not to retract their paper.} set out to find explicit solutions to a stochastic Burgers' equation and compare these with numerical simulations. If successful, this would be a welcome contribution to the field. In order to solve the equation, the authors use Ansatz \textbf{A}. This Ansatz goes against our intuition of how a stochastic travelling wave should look and indeed, the computations will show that this Ansatz is invalid. However, the authors ignore this and go on to compute a range of solutions, mostly unphysical, and all of them are not solutions of the stochastic Burgers' equation. The paper contains figures that should validate the solutions by comparison with numerical solutions, but it is not clear what these figures mean; they definitely do not show a comparison of an analytical and numerically computed travelling wave. Furthermore, the authors argue that their Euler-Maruyama scheme is unconditionally stable, which it is not~\cite{lord2014book}.\\
\begin{remark}
A warning here is in place. I use the notation $\beta_t$ to indicate Brownian motion, which is often written as $W_t$ or $B_t$ in the literature. However, in many of the papers discussed here, the subscript $t$ is used to denote the time derivative. Hence, in \cite{baber2024exact}, $\beta(t)$ is Brownian motion and $\beta_t$ is white noise. I will not follow this notation here. As an example of how this leads to confusion, see the PubPeer discussion for \cite{baber2024exact}. 
\end{remark}
\noindent\textbf{Travelling Wave Ansatz}\\
The stochastic Burgers' equation as studied in \cite{baber2024exact} is given by
\begin{align}
\label{eq:SBE}
d\U=[\nu\U_{xx}-\U\U_x]dt+\sigma \U d\beta_t.
\end{align}
For the deterministic case, i.e. $\sigma=0$, the travelling wave $(\Phi,c)$ can be found by switching to the coordinate $\xi=x-ct$, resulting in the Travelling Wave Equation
\begin{align}
\label{eq:TWE}
\nu\Phi''+c\Phi'-\Phi\Phi'=0.
\end{align}
Note that each constant function is a solution. Upon fixing the boundary conditions $\Phi(-\infty)=f^+$ and $\Phi(\infty)=f^-$ with $f^+>f^-$, the travelling wave is given by
\begin{align}
\label{eq:SolTWE}
    \Phi(\xi)= c - \frac{f^+-f^-}{2}\tanh\left(\frac{f^+-f^-}{4\nu}\xi\right), \hspace{5mm} c=\frac{f^++f^-}{2}.
\end{align}
This solution can be found by integrating once and using partial fraction decomposition. Note here that boundary conditions are never discussed in \cite{baber2024exact}, but they determine the wave speed. To find a solitary wave, as the title of the paper claims, we must have that $f^+=f^-$, but this results in the trivial zero solution. As far as I know, there are no solitary travelling waves in Burgers' equation, and in \cite{baber2024exact} none are shown, despite the title. Furthermore, note that in \cite{baber2024exact} the travelling wave coordinate is defined as $\xi=lx-ct$, but this $l$, which they call the amplitude of the wave, remains a free parameter throughout the paper. Also note that the authors refer to $c$ not as the travelling wave speed, but as the speed of light. Lastly, when $f^-=0$, the solution in \eqref{eq:SolTWE} can also be written as
\begin{align}
    \Phi(\xi)= \frac{2 c}{1+e^{cx/\nu}}. 
\end{align}
Hence, solutions like $\U_5$, $\U_{15}$ and $\U_{16}$ are just another way of writing \eqref{eq:SolTWE}.\\
The authors now assume that a solution to \eqref{eq:SBE} can be written as
\begin{align}
\U(x,t)=\Theta(\rho)e^{\sigma\beta_t-\frac{\sigma^2}{2}t}, \hspace{1cm}\rho=lx-ct.
\end{align}
Here, $\Theta$ is assumed to be a deterministic function.
Following the computations in \cite{baber2024exact} using It\^o calculus, this results in the equation
\begin{align}
\label{eq:StochTW}
-c\Theta'+l\Theta\Theta'e^{\sigma\beta_t-\frac{\sigma^2}{2}t}-\nu l^2\Theta''=0.
\end{align}
In order to satisfy the Ansatz, we must show that a deterministic solution exists to the equation above. As far as I know, these do not exist apart from trivial constant solutions. Therefore, this Ansatz recovers only spatially homogeneous solutions. 

However, the authors claim that deterministic solutions to \eqref{eq:StochTW} can be found by replacing $e^{\sigma\beta_t-\frac{\sigma^2}{2}t}$ with the average $E[e^{\sigma\beta_t-\sigma^2t}]=1$. When we set the undetermined variable $l=1$,  Eq.~\eqref{eq:StochTW} hence reduces to the classic TWE \eqref{eq:TWE}. Summarising, the authors claim that any solution to \eqref{eq:TWE} is also a solution to the stochastic Burgers' equation when we multiply it by $e^{\sigma\beta_t-\frac{\sigma^2}{2} t}$. As this approach does not depend on the specific equation, we hence conclude that for each reaction-diffusion equation that has a travelling wave, the stochastic version with $\sigma\U d\beta_t$ added can be solved by multiplying the travelling wave with geometric Brownian motion. This contradicts all available literature on stochastic travelling waves. Of course, we are not surprised, as the deterministic solutions of \eqref{eq:TWE} are not solutions of \eqref{eq:StochTW}, which can be verified by direct substitution.  

Here are some more examples that use Ansatz \textbf{A}, combined with the technique of replacing $e^{\sigma\beta_t-\frac{\sigma^2}{2}t}$ by 1: \cite{al2023exact, mohammed2021exact, kaya2025analytical, baber2025solitary}.

\section{Several classic PDEs}
In this section, we will discuss the following question: Are there stochastic extensions of classical PDEs that can be solved by multiplying the deterministic solution by a stochastic exponent? To simplify this rather broad question, we restrict ourselves to parabolic equations forced by a 1D Wiener process. 

\subsection{Burgers' Equation}

To stay close to the case studied above, we start with a stochastic Burgers' equation:
\begin{align}
\label{eq:SBEGen}
d\U=[\nu\U_{xx}-\U\U_x]dt+\sigma g(\U) d\beta_t.
\end{align}
This equation is more general than \ref{eq:SBE}, as $\U$ has been replaced by $g(\U)$. Generalizing Ansatz \textbf{A}, we propose 
\begin{align}
\label{eq:StochBurg}
    \Upsilon(x,t)=\Theta(x-ct)h(\beta_t,t),
\end{align}
for some, and this is essential, deterministic function $\Theta$.
Using It\^o calculus, we find
\begin{align}
\label{eq:ItoCalcBurgers}
    d\Upsilon =[-c\Theta'(x-ct)h(\beta_t,t)+\Theta(x-ct)h_2(\beta_t,t)+\frac{1}{2}h_{11}(\beta_t,t)]dt+\Theta(x-ct)h_1(\beta_t,t)d\beta_t.
\end{align}
Comparing the stochastic terms then results in
\begin{align}
    \sigma g(\Theta(x-ct)h(\beta_t,t))=\Theta(x-ct)h_1(\beta_t,t).
\end{align}
As this must hold for all $x$, we find that $g(z)=z$ and hence the equation above reduces to
\begin{align}
    \sigma h(\beta_t,t)=h_1(\beta_t,t).
\end{align}
Solving this equation results in 
\begin{align}
    h(\beta_t,t)=e^{\sigma \beta_t}\tilde h(t).
\end{align}
Comparing the deterministic terms in \sref{eq:SBEGen} and \sref{eq:ItoCalcBurgers} now results in 
\begin{align*}
    \nu\Theta''h(\beta_t,t)-\Theta\Theta'h^2(\beta_t,t)=-c\Theta'(x-ct)h(\beta_t,t)+\Theta(x-ct)h_2(\beta_t,t)+\frac{\sigma^2}{2}\Theta h(\beta_t,t).
\end{align*}
Hence,
\begin{align}
\label{eq:TheProblem}
    \nu\Theta''-\Theta\Theta'h(\beta_t,t)=-c\Theta'+\Theta\tilde h'(t)/\tilde h(t)+\frac{\sigma^2}{2}\Theta.
\end{align}
The only choice to make this equation deterministic is to choose $\Theta'=0$.  

We conclude that this Ansatz only results in trivial solutions. Zooming out, what prevents this Ansatz from working? For this Ansatz to be effective, it is essential that any dependence on $t$ or $\beta_t$ disappears from the TWE. It is clear from Eq. \sref{eq:TheProblem}, that this only works for linear equations. However, as we are studying nonlinear waves, this is impossible. This computation is not limited to Burgers' equation. For equations like Nagumo, FitzHugh-Nagumo, Korteweg-De Vries, Camassa-Holm etc., Ansatz \textbf{A} will always result in an equation like Eq.~\sref{eq:TheProblem}, which cannot be solved by nontrivial deterministic solutions.

\subsection{NLS equation}
Is a successful application of Ansatz \textbf{A} then completely impossible? No, imaginary numbers can save us.
We start with stating the classic undimensionalised focusing cubic NLS equation:
\begin{align}
    i\psi_t=-\psi_{xx}-|\psi|^2\psi.
\end{align}
We propose a classic Ansatz of the form $\psi(x,t)=\Phi(x-ct)e^{ikx-i\omega t}$. Plugging this back into the NLS equation then results in
\begin{align}
    i(-c\Phi'-i\omega \Phi)e^{ikx-i\omega t}=-(\Phi''+2ik\Phi'-k^2\Phi)e^{ikx-i\omega t}-\Phi^3e^{ikx-i\omega t}
\end{align}
Equating the real and imaginary parts then gives us
\begin{align}
\label{eq:NLSTWE}
    \begin{split}
        \omega\Phi&=-\Phi''+k^2\Phi-\Phi^3,\\
        -c\Phi'&=-2k\Phi'.
    \end{split}
\end{align}
Hence, $c=2k$. When we multiply the equation once by $\Phi'$ and integrate, we find
\begin{align}
       \frac{\omega}{2}\Phi^2&=-\frac{1}{2}\Phi'^2+\frac{1}{2}k^2\Phi^2-\frac{1}{4}\Phi^4+C.
\end{align}
Note that $C=0$ as $\Phi(+\infty)=0$, hence
\begin{align}
       \Phi'^2-(k^2-\omega)\Phi^2+\frac{1}{2}\Phi^4=0.
\end{align}
Now, $\Phi=A\text{sech}(B\xi)$ solves this equation when $B=A/\sqrt{2}$ and $k^2-\omega=A^2/2$. Hence, we have a 2D family of solutions, and 3D if we include translation invariance. Why does this approach work? Because the complex exponent disappears from the equation of the soliton. Hence, for every nonlinear term of the form $|u|^{2\alpha}u$, this Ansatz works.

We now extend the NLS equation to a stochastic version:
\begin{align}
\label{eq:StochNLS}
    id\psi=[-\psi_{xx}-|\psi|^2\psi]dt+\sigma \psi d\beta_t
\end{align}
Such an equation (with infinite-dimensional noise) has already been studied numerically before \cite{debussche2002numerical}. 
Following \cite{ali2024}, we propose the following Ansatz, which we will refer to as Ansatz \textbf{B}:
\begin{align}
    \psi=\Phi(x-ct)e^{i(kx-\omega t)+\sigma\beta_t-\frac{\sigma^2}{2}t}
\end{align}
There are two reasons Ansatz \textbf{B} cannot work:
\begin{enumerate}
    \item Using It\^o calculus, we find that the stochastic part of $id\psi$ is $i\psi d\beta_t$, not $\psi d\beta_t$. 
    \item To find an equation like \sref{eq:NLSTWE}, it is essential that the exponent in $|\psi|^2\psi$ cancels against the exponent of the linear terms. Hence, the exponent must be purely imaginary. 
\end{enumerate}
We could circumvent problem 1 by assuming a complex noise term in Eq.~\sref{eq:StochNLS}, but that would still leave problem 2. 
Therefore, we update our Ansatz to
\begin{align}
    \psi=\Phi(x-ct)e^{i(kx-\omega t-\sigma\beta_t+\frac{\sigma^2}{2}t)},
\end{align}
which we will refer to as Ansatz \textbf{C}. 
This then results in 
\begin{align}
    d\psi=[-c\Phi'+(-i\omega +i\frac{\sigma^2}{2})\Phi+\frac{(-i\sigma)^2}{2}\Phi]e^{i(kx-\omega t-\sigma\beta_t+\frac{\sigma^2}{2}t)}dt-i\sigma\psi d\beta_t.
\end{align}
The stochastic parts of the NLS equation and $id\psi$ now match. For the deterministic part, we find
\begin{align}
    i(-c\Phi'-i\omega \Phi-\frac{\sigma^2}{2}\Phi+i\frac{\sigma^2}{2}\Phi)=-(\Phi''+2ik\Phi'-k^2\Phi)-\Phi^3.
\end{align}
Splitting into real and imaginary parts, we get
\begin{align}
\label{eq:NLSunintegrable}
    \begin{split}
        \omega \Phi-\frac{\sigma^2}{2}\Phi&=-\Phi''+k^2\Phi-\Phi^3\\
        -c\Phi'-\frac{\sigma^2}{2}\Phi&=-2k\Phi'
    \end{split}
\end{align}
The second equation could not possibly have a nontrivial integrable solution; hence, we find no solitons. 

Stepping back, could we have known from the beginning that this approach would not work? Yes, given the numerical results in \cite{debussche2002numerical}, we know that the amplitude decays over time. Hence, we expect that this time-dependence must be explicitly modelled by the Ansatz, which is not the case for Ansatz \textbf{C}.  

Can we fix the problem in Eq.~\sref{eq:NLSunintegrable} above? First, we note that we do not need the It\^o correction in the Ansatz, as it ends up in the real part while the It\^o correction ends up in the imaginary part. However, we could balance it by ad-hoc introducing a deterministic term $i\frac{\sigma^2}{2}\psi$ into the NLS equation, which the reader will recognise as an It\^o-Stratonovich correction term.

\subsubsection{Stratonovich interpretation}
It seems that the troubles come from the It\^o calculus, so let us study the following equation
\begin{align}
\label{eq:NLSstrat}
    id\psi=[-\psi_{xx}-|\psi|^2\psi]dt+\sigma \psi\circ d\beta_t,
\end{align}
and introduce a modified Ansatz \textbf{Cs}:
\begin{align}
    \psi(x,t)=\Phi(x-ct)e^{i(kx-\omega t-\sigma\beta_t)}=\psi_\mathrm{det}(x,t)e^{-i\sigma\beta_t}.
\end{align}
Applying Stratonovich calculus then gives us 
\begin{align}
\label{eq:NLSstratcalc}
    d\psi=[-c\Phi'-i\omega\Phi]e^{i(kx-\omega t-\sigma\beta_t)}dt-i\sigma\psi\circ d\beta_t.
\end{align}
Equating \sref{eq:NLSstrat} and \sref{eq:NLSstratcalc}, we get
\begin{align}
    \begin{split}
        \omega \Phi&=-\Phi''+k^2\Phi-\Phi^3,\\
        -c\Phi'&=-2k\Phi',
    \end{split}
\end{align}
which is identical to the deterministic equation. Hence, Ansatz \textbf{Cs} indeed works for the NLS equation in the Stratonovich interpretation.  
Is this a general lesson? Can we find solutions to other SPDEs in the Stratonovich interpretation? No, because even in the Stratonovich interpretation, it is still essential that the equations above are linear in the stochastic exponent.  
An important remark here is that in many cases in the literature, it is not clear whether the It\^o or the Stratonovich interpretation is used. Or when it is mentioned, it might not be actually used, for example, in \cite{alkhidhr2023new,arnous2024investigating}, it is stated that the It\^o interpretation is used, but Ansatz \textbf{Cs} is used for NLS-type equations, and the computations are done using Stratonovich calculus.

\section{Where to next}
It is an intriguing question why there are so many of these articles around, given the obvious problems. Answering this question would lead to speculation on the intent of the authors; are they honestly mistaken or fraudulent? I do not wish to discuss the motivations of individual researchers here, especially as it can be almost too easy to criticise others from a cosy, western university. We can, and should, however, criticise the journals that publish this kind of work, especially Q1 journals. Unfortunately, sending editors of journals a request for retraction does not work, it only results in `promises to investigate', but none have materialised so far.

I list here some responses from editors that are from cosy, western universities:
\begin{itemize}
    \item Dr. Rafal Marszalek, editor for Scientific Reports, promised an investigation into \cite{baber2024exact}, but after two years and several repeated promises, there are no updates. I also submitted a `Matters Arising' to Scientific Reports about this article, but this was declined after half a year of waiting. According to the journal, the problems with the paper are not suitable for a Matters Arising, and again an investigation was promised.
    \item Dr. Gary Ren, editor for Frontiers in Physics, did initially not respond to questions on \cite{alkhidhr2023new}. After I sent him the first version of this document, he stated that ``It would not be appropriate for me, as an individual editor, to adjudicate these matters or speculate on referees’ intent outside those formal processes." This is of course nonsense, as adjudicating these matters is exactly the job description of an editor. 
    \item AIMS mathematics stated that there were no reasons to retract \cite{trouba2025soliton} after consultation with Prof. Brzezniak. Prof. Brzezniak has not replied yet to inquiries. 
    \item Dr. Helmut Abels, editor for ZAMM, did not react to questions immediately after publication of \cite{al2023exact}, nor to a draft of this note. I asked for the input of Guido Schneider, also editor for ZAMM. He did recognise this paper as ``complete nonsense".  Now, Helmut Abels did respond, stating that he would be willing to publish a comment by me alongside the paper. So I did, but the computer system and the editor in chief had some difficulty in processing a `comment'. Presently, the comment is under review. \item For PLOS One, I tried a different approach. I flagged \cite{ali2024} as plagiarism, as it is almost identical to \cite{mohammed2021exactSchrodinger}. The journal promised to investigate the paper. 
    \item Wen-Xiu Ma, professor at the University of South Florida, does not believe he has to answer questions on papers like \cite{rasid2025abundant} because he is not the corresponding author, even though a) that argument does not make sense, b) his email is on the paper and no explicit corresponding author is mentioned.  
\end{itemize}

It is not practical to write to every individual editor explaining the issues with each individual paper. The point of this document is hence to reduce the effort for other scientists to report these types of errors, either on a site like Pubpeer or directly to the editor. 

Ideally, this document becomes interactive, in the sense that I invite other researchers to contribute, e.g. by describing similar Ans\"atze in other equations.\\

\noindent\textbf{Competing Interests}\\
There are no competing interests. \\

\noindent\textbf{Acknowledgements}\\
The work of CH is funded by the Dutch Institute for Emergent Phenomena (DIEP) at the University of Amsterdam via the program Foundations and Applications of Emergence (FAEME).

\bibliographystyle{klunumHJ}
\bibliography{ref}

\end{document}